\documentclass[11pt,oneside,reqno]{amsart}
\usepackage{amssymb}
\usepackage{amsmath}
\usepackage{amscd}
\usepackage{amsfonts}
\usepackage{enumitem}
\usepackage{mathrsfs}
\usepackage{stmaryrd}
\usepackage{tikz}
\usepackage{hyperref}

\usetikzlibrary{arrows}
\usepackage[T1]{fontenc}
\allowdisplaybreaks

\makeatletter
\DeclareFontFamily{U}{tipa}{}
\DeclareFontShape{U}{tipa}{m}{n}{<->tipa10}{}
\newcommand{\arc@char}{{\usefont{U}{tipa}{m}{n}\symbol{62}}}%

\newcommand{\arc}[1]{\mathpalette\arc@arc{#1}}

\newcommand{\arc@arc}[2]{%
  \sbox0{$\m@th#1#2$}%
  \vbox{
    \hbox{\resizebox{\wd0}{\height}{\arc@char}}
    \nointerlineskip
    \box0
  }%
}
\makeatother

\usepackage{bbm}

\theoremstyle{definition}
\newtheorem{theorem}{Theorem}[section]

\newtheorem{thm}[theorem]{Theorem}

\newtheorem{prop-def}{Proposition-Definition}[section]

\newtheorem{rema}[theorem]{Remark}

\newcommand{\N}{{\mathbb N}}
\newcommand{\C}{{\mathbb C}}
\newcommand{\Z}{{\mathbb Z}}

\newcommand{\Q}{{\mathbb Q}}
\newcommand{\g}{{\mathfrak g}}

\newcommand{\one}{\mathbf{1}}
\renewcommand{\d}{\mathbf{d}}
\newcommand{\wt}{\mbox{\rm wt}\ }

\begin{document}

\setlength{\oddsidemargin}{0cm} \setlength{\evensidemargin}{0cm}
\baselineskip=18pt

\title{Deformation rigidity of some simple affine VOAs}

\author{Andrew R. Linshaw and Fei Qi}

\begin{abstract}
    In this paper, we prove that simple affine vertex operator algebras with positive integer levels admit only trivial first-order deformations. Therefore, the deformation rigidity conjecture of strongly rational vertex operator algebras holds for these cases. We also show that the same holds for simple affine vertex operator algebra of $\mathfrak{sl}_2$ at the non-integral admissible level $-4/3$. Therefore, neither $C_2$-cofiniteness nor rationality is a necessary condition for deformation rigidity of VOAs. We conjecture that the same should hold for every simple affine VOA that does not coincide with the corresponding universal affine VOA. 
\end{abstract}

\maketitle

\section{Introduction}

The moduli space of conformal field theories is very important in both physics and mathematics and connects to a number of branches of mathematics. But we know very little about this moduli space, not even the correct topology. It has been a common belief among physicists since the 1990s that among all two-dimensional conformal field theories, the strongly rational ones are isolated. The belief is formulated into a mathematical conjecture by Yi-Zhi Huang in \cite{H-open}, Problems 4.1 and 4.2. It is therefore important to study the deformations of vertex algebras. One of the conjectures in this study is what we call the deformation rigidity conjecture: the deformation of a rational vertex operator algebra that preserves the grading and the vacuum can only be trivial.

Classifying the first-order deformations has been a long-standing challenge, mainly due to technical convergence conditions. The first significant result did not appear until 2024, when the second author, in collaboration with Vladimir Kovalchuk, overcame the difficulty in convergence, and gave a programmable algorithm for classifying the first-order deformations for every freely generated vertex algebra (see \cite{KQ-1st-def}). 

One of the key theorems in \cite{KQ-1st-def} states that for the vertex operator $Y_1$ associated with the first-order deformation, if for any generators $a, b$ of the vertex algebra, $Y_1(a,x)b$ contains no negative powers, then $Y_1$ gives a trivial deformation. This theorem allows us to focus only on the singular part of $Y_1(a, x)b$ with only finitely many structural constants. \cite{KQ-1st-def} proved the key theorem by constructing the linear map corresponding to the trivial deformation based on the PBW-type basis on a freely generated vertex algebra. The proof fails if there are relations among the basis vectors. 

In the case of simple affine vertex operator algebras $L_k(\g)$ with the level $k$ being a positive integer level, we show that the same result holds by using a theorem of Bong Lian \cite{Lian} to circumvent the direct construction in \cite{KQ-1st-def}. Then, by combining some of the computations in \cite{KQ-1st-def}, we finish the proof of the deformation rigidity conjecture for these particular examples of rational vertex operator algebras.

We also prove the same result for $L_{-4/3}(\mathfrak{sl}_2)$, a particular case of a simple affine VOA at a non-integral admissible level. Note that $L_{-4/3}(\mathfrak{sl}_2)$ is neither $C_2$-cofinite nor rational, so neither of these conditions is necessary for deformation rigidity. We conjecture that the deformation rigidity should hold for all simple affine VOAs $L_k(\g)$ that do not coincide with the corresponding universal affine VOA $V^k(\g)$. If the conjecture holds, then together with the results in \cite{KQ-1st-def}, we have 
$$\dim H^2_{1/2}(L_k(\g), L_k(\g)) = \begin{cases}
    1 & \text{ if }L_k(\g) = V^k(\g)\\
    0 & \text{ if }L_k(\g) \neq V^k(\g)
\end{cases}.$$
The proof, however, might require a completely different method. 

The paper is organized as follows. Section 2 reviews the notion of first-order deformations and their equivalence classes. Section 3 shows that a regular $Y_1$ is trivial. In Section 4, we give a proof that $Y_1$ is regular for the simple affine VOA $L_k(\g)$ with $k\in \Z_+$. In Section 5, we give a proof that $Y_1$ is regular for the simple affine VOA $L_{-4/3}(\mathfrak{sl}_2)$. Both Section 4 and 5 substantially depend on the formulas and computations in Section 4.1 and 6.2 of \cite{KQ-1st-def}. 

\textbf{Acknowledgement: } The second author would like to thank Yi-Zhi Huang for his long-term support. The first author is supported by NSF Standard Grant DMS-2401382 and Simons Foundation Grant MPS-TSM-00007694. Both authors would like to thank Yi-Zhi Huang, Shashank Kanade, and Paul Johnson for their comments on the early versions of the draft. 

\section{Preliminaries}

\subsection{First-order deformations}\label{first-order-deform-subsec}
Let $(V, Y, \one)$ be a grading-restricted vertex algebra in the sense of \cite{H-Coh}. This means a $\Z$-graded vertex algebra with usual grading restrictions, and with only the grading operator $\d=L(0)$ and $D = L(-1)$. For $u, v\in V$, we keep Borcherds' notation $u_n$ for the modes of $Y(u, x)$, i.e., 
$$Y(u, x)v = \sum_{n\in \Z} u_n v x^{-n-1}$$
A first-order deformation of $V$ is defined by a vertex operator 
\begin{align*}
    Y_1: \ & V\otimes V \to V((x))\\
    & u\otimes v \mapsto Y_1(u, x)v = \sum_{n\in \Z} u_n^{def}v x^{-n-1}. 
\end{align*}
such that the vector space $V^t = \C[t]/(t^2)\otimes_\C V$, with the vertex operator
$$Y^t(u, x) v = Y(u, x)v + t Y_1(u, x)v,$$
and the vacuum element $\one$, forms a grading-restricted vertex algebra over the base ring $\C[t]/(t^2)$ (see \cite{H-1st-2nd-Coh} and \cite{KQ-1st-def} for more details). Equivalently, $Y_1$ should satisfy the following conditions: 
\begin{enumerate}[leftmargin=*]
\item For every $u, v\in V$, the mode $u^{def}_n$ satisfies
$$\wt u^{def}_n v = \wt u - n - 1 + \wt v. $$
\item For every $v\in V$, 
\begin{align}
    Y_1(\one, x) v = 0, Y_1(u, x) \one = 0. \label{Y-1-on-vac}
\end{align}
In particular, the $D$-operator 
$$D (u+tv) = \lim_{x\to 0}\frac{d}{dx} Y^t(u+tv, x) \one = \lim_{x\to 0}\frac{d}{dx} Y(u+tv, x)\one = u_{-2}\one + t v_{-2}\one $$on $V^t$ is defined purely with modes of $Y$, regardless of $Y_1$. 
\item For every $u, v\in V$, the $D$-derivative-commutator formula holds,  
$$[D, Y_1(u, x)] v = Y_1(Du, x) v = \frac{d}{dx} Y_1(u, x) v. $$
\item For $u, v\in V$, the following Jacobi-like identity holds,
\begin{align*}
    & x_0^{-1}\delta\left(\frac{x_1-x_2}{x_0}\right) \left(Y_1(u, x_1)Y(v, x_2) + Y(u, x_1)Y_1(v, x_2)\right)\\
    & - x_0^{-1}\delta\left(\frac{-x_2+x_1}{x_0}\right)\left(Y_1(v, x_2)Y(u, x_1) + Y(v, x_2)Y_1(u, x_1)\right)\\
    = \ & x_2^{-1}\delta\left(\frac{x_1-x_0}{x_2}\right)\left(Y_1(Y(u,x_0)v, x_2) + Y(Y_1(u, x_0)v, x_2)\right).
\end{align*}
\end{enumerate}

\begin{rema}
For $u, v\in V$ and $m, n\in \Z$, by taking an appropriate residue of the Jacobi-like identity, we have the following commutator condition.
\begin{align}
    [u_m^{def}, v_n] + [u_m, v_n^{def}] = \sum_{\alpha\geq 0}\binom{m}{\alpha} \left((u_\alpha^{def} v)_{m+n-\alpha} + (u_\alpha v)^{def}_{m+n-\alpha}\right).\label{Comm-reln}
\end{align}
The commutator relation (\ref{Comm-reln}) should hold for every $u, v\in V$ and $m, n\in \Z$. In case $V$ is a freely generated vertex algebra, \cite{KQ-1st-def} proves, via a lengthy process, that it suffices for (\ref{Comm-reln}) to hold for every pair of generators $u, v$ and every $m, n\in \N$. The result does not naively generalize to non-freely generated vertex algebras. But we would not need it if no nontrivial first-order deformations exist.

\end{rema}
\subsection{Equivalent first-order deformations} Two first-order deformations $(V^t, Y^{t,(1)})$ and $(V^t, Y^{t,(2)})$ are equivalent, if there exists a $\C[t]/(t^2)$-linear vertex algebraic isomorphism $f^t: V \oplus tV = V^t \to V^t =V\oplus tV $ whose restriction on $V$ is of the form
$$f^t|_V = 1_V + t f_1, $$
where $f_1:V\to V$ is a $\C$-linear grading-preserving map. In other words, for $u,v\in V$, 
$$f^t(u+tv) = u + t(v+f_1(u)). $$

In the papers \cite{H-Coh} and \cite{H-1st-2nd-Coh}, Huang showed that the equivalence classes of first-order deformations are described by the second cohomology $H_{1/2}^2(V, V)$ of vertex algebras. Here we shall not give a detailed review of the definitions of 2-cocycles and 2-coboundaries, but only mention one consequence. Fix an $Y_1: V \otimes V \to V((x))$. If there exists a grading preserving map $\phi: V \to V$ such that $\phi(\one)=0$ (see \cite{H-1st-2nd-Coh-Add}), and for every $u, v\in V$,
$$Y_1(u, x)v = Y(\phi(u), x) v - \phi(Y(u,x)v)+ Y(u, x)\phi(v),$$
then the first-order deformation given by $Y_1$ is equivalent to that given by $0$. In this case, we say $Y_1$ is a trivial first-order deformation. 

\subsection{Simple affine VOA}
In this paper, we focus on the particular example where $V$ is the simple affine VOA (see \cite{Frenkel-Zhu}, \cite{LL} for more details). Let $\g$ be a finite-dimensional simple Lie algebra with root system $\Phi$. Let $\theta$ be the highest root of $\g$. For every $\alpha\in \Phi$, let $e_\alpha, f_\alpha, h_\alpha$ be the basis vectors of the $\mathfrak{sl}_2$-subalgebra of $\g$ associated with $\alpha$. When $\alpha = \theta$, we abbreviate $e_\theta, f_\theta, h_\theta$ as $e, f, h$. For a fixed real number $k$, let $V^k(\g)$ be the universal affine VOA associated with $\g$ with level $k$. Let $L_k(\g)$ be the simple quotient of $V^k(\g)$. In case $k+h^\vee\in \Q_{\geq 0}\setminus\{1/m: m\in \Z_+\}$, $L_k(\g)$ is a nontrivial quotient of $V^k(\g)$ (see \cite{Gorelik-Kac} for details). In case $k$ is admissible, $L_k(\g)$ is the quotient of $V^k(\g)$ by an ideal generated by singular vectors. (see \cite{Kac-Wakimoto}). In case $k$ is a positive integer,  $L_k(\g)$ is the quotient of $V^k(\g)$ by the ideal generated by $e(-1)^{k+1}\one$ (see \cite{LL}). 

Hereafter, we set $V = L_k(\g)$, $V^t = \C[t]/(t^2)\otimes_\C V = V \oplus t V$ the first-order deformation of $V$ with the vertex operator $Y^t = Y + t Y_1$. For $a\in \g$, we use $a(m)$ and $a^{def}(m)$ to denote the modes in the vertex operators $Y(a,x)$ and $Y_1(a,x)$, respectively.

\section{Triviality of a regular $Y_1$}

For a freely generated vertex algebra, \cite{KQ-1st-def} shows that if the $Y_1$-operator on every pair of generators is regular, then $Y_1$ defines a trivial first-order deformation. Here, using a result in \cite{Lian}, we generalize it to the simple affine VOAs that are not freely generated. 

\begin{thm}\label{thm-3-1}
    Let $V$ be a simple affine VOA. If for every $a, b\in \g$, $Y_1(a(-1)\one, x)b(-1)\one$ is regular, i.e., it contains no negative powers, then $Y_1$ defines a trivial first-order deformation. 
\end{thm}

\begin{proof}
    Consider the $\C$-subalgebra $W$ of $V^t$ generated by the vacuum $\one$ and the weight-1 subspace $V_{(1)} \subset V_{(1)} \oplus t V_{(1)} \subset V^t$. Since $Y_1$ is regular, for $a, b\in V_{(1)}$, the singular part of $Y^t(a,x)b$ coincides with that of $Y(a,x)b$. From Theorem 4.11 of \cite{Lian}, $W$ is isomorphic to a quotient of affine VOA. Clearly, $W_{(1)}=V_{(1)}$ is isomorphic to $\g$. So $W$ is isomorphic to a quotient of $V^k(\g)$. In particular, $W_{(0)} = \C \one$. 

    Consider the projection $\pi: V^t \to V$ modulo $t$, i.e., $\pi(v + tv^{(1)}) = v$. Clearly, $\pi$ is a vertex algebra homomorphism. Let $\pi|_W$ be its restriction on $W$. Since $W_{(1)} = V_{(1)}, \pi|_W$ is surjective. Let $K = \ker \pi|_W$. Clearly, $K$ is an ideal of $W$ contained in $tV \cap W$. Let 
    $$I = \{v\in V: tv \in K\}.$$
    We show that $I$ is an ideal of $V$. Fix any $u\in V$ and $v\in I$. Let $w = u + t w^{(1)}\in W$ so that $\pi(w) = u$. Since $tv$ is an element of the ideal $K\subset W$, for every $n\in \Z$, $(w_n + tw_n^{def})\cdot tv\in K$, i.e., 
    $$(w_n + t w_n^{def}) \cdot tv = (u_n + t w^{(1)}_n + t w_n^{def}) \cdot tv = t u_n v$$
    is an element in $K$. Thus, $u_n v\in I$. 

    Since $V$ is simple, either $I = 0$ or $I = V$. If $I = V$, then in particular, $\one\in I$, i.e., $t\one \in W$, contradicting $W_{(0)} = \C \one$. Therefore, $I = 0$. Consequently, $K = 0$. So $\pi|_W$ is an isomorphism. We thus proved that $W$ is isomorphic to $V = L_k(\g)$. 

    Let $f = (\pi|_W)^{-1}: V \to W$ be the isomorphism. From $\pi \circ f = 1_V$, we know that for every $v\in V$
    $$f(v) = v - t \phi(v)$$
    where $\phi: V \to V$ is a homogeneous linear map. Fix any $u, v\in V$. We read the coefficient of $t$ from 
    $$f(Y(u,x)v) = Y^t(f(u), x)f(v),$$
    to obtain 
    $$-\phi(Y(u,x)v) = Y_1(u, x) v - Y(\phi(u), x)v - Y(u, x) \phi(v). $$
    Therefore, $Y_1(u,x)v$ coincides with the trivial deformation given by $\phi$, thus trivial. 
\end{proof}

\section{The positive integer level case}\label{sec-int-level}

In this section we show that $Y_1$ is regular. The computation in Section 6.2 of \cite{KQ-1st-def} generalizes to the current situation, showing that the first-order pole of $Y_1$ is represented by a trivial first-order deformation, and for every $a, b\in \g$, 
\begin{align}
    Y_1(a(-1)\one, x)b(-1)\one = c\cdot \langle a, b\rangle \one x^{-2}+ \sum_{m\geq 1} a^{def}(-m)b(-1)\one x^{m-1}\label{Y_1-modes}
\end{align}
for some constant $c\in \C$ (which is denoted by $B_\g$ in \cite{KQ-1st-def}). Here $\langle \cdot, \cdot \rangle$ is the invariant bilinear form on $\g$. So, it suffices to show that $c=0$.

The key idea is to consider the coefficient of $e(-1)^k\one$ in the expression of $f^{def}(1)e(-1)^{k+1} \one$. This expression is zero because $e(-1)^{k+1} \one = 0$ in $V$. On the other hand, we may use the commutator condition (\ref{Comm-reln}) to calculate it. This shall establish an equation of $c$. 




\subsection{The mode $f^{def}(1)$}
We now consider the commutator condition (\ref{Comm-reln}) with $u= f(-1)\one, v= e(-1)\one$, $m=1$, and $n = -1$. Then the commutator condition reads 
\begin{align*}
    & [f^{def}(1), e(-1)] + [f(1), e^{def}(-1)] \\
    = \ & \left(f^{def}(0) e(-1)\one\right)_{1-1-0} +\left( f(0)e(-1)\one\right)^{def}_{1-1-0} + \left(f^{def}(1) e(-1)\one\right)_{1-1-1} +\left( f(1)e(-1)\one\right)_{1-1-1}^{def}\\
    = \  & 0 - h^{def}(0)+ (c \cdot \one)_{-1} + (k \cdot \one)_{-1}^{def} = -h^{def}(0) + c. 
\end{align*}
Unraveling the left-hand-side, we have 
\begin{align}
    f^{def}(1) e(-1) - e(-1) f^{def}(1) + f(1) e^{def}(-1) - e^{def}(-1)f(1) = -h^{def}(0) + c. \label{mode-1}
\end{align}
%
Let $i\in \Z_+, 1\leq i \leq k+1$. Apply (\ref{mode-1}) to $e(-1)^{i-1}\one$, we see that
\begin{align}
    f^{def}(1)e(-1)^i\one = \ & e(-1)f^{def}(1)e(-1)^{i-1}\one - f(1) e^{def}(-1)e(-1)^{i-1}\one + e^{def}(-1)f(1)e(-1)^{i-1}\one\nonumber\\
    & - h^{def}(0)e(-1)^{i-1}\one + c\cdot e(-1)^{i-1}\one. \label{Recursion}
\end{align}
In particular, if $i = k+1$, then left-hand-side is zero. 

\subsection{The term with $h^{def}(0)$}\label{h0-def-sec} We start by analyzing $h^{def}(0) e(-1)^{i}\one$ and show that it is zero for every $1\leq i \leq k$. 

We argue by induction. From (\ref{Y_1-modes}), we see that 
    $$h^{def}(0) e(-1)\one = 0.$$
    Thus, the conclusion holds for $i=1$. 
    
    Assume the conclusion holds for smaller $i$. In the commutator relation (\ref{Comm-reln}), we set $u = h(-1)\one, v = e(-1)\one, m = 0$, and $n = -1$. Using (\ref{Y_1-modes}), 
    \begin{align}
        & h^{def}(0) e(-1) - e(-1)h^{def}(0) + h(0)e^{def}(-1) - e^{def}(-1)h(0) \nonumber\\
        = \ &  (h^{def}(0) e(-1)\one)_{0-1-0} + (h(0) e(-1)\one)^{def}_{0-1-0} =  2 e^{def}(-1). \label{h0-def-comm-condn}
    \end{align}
    Apply these modes on $e(-1)^{i-1}\one$, we see that
    \begin{align}
        & h^{def}(0) e(-1)^i\one - e(-1)h^{def}(0) e(-1)^{i-1}\one + h(0)e^{def}(-1)e(-1)^{i-1}\one  - e^{def}(-1)h(0)e(-1)^{i-1}\one  \nonumber\\
        = \ &  2 e^{def}(-1)e(-1)^{i-1}\one. \label{h0-def-e-i-1}
    \end{align}
    Since there are no relations in conformal weights at most $k$, the coboundary construction Section 3 of \cite{KQ-1st-def} applies verbatim in this range. Replacing $Y_1$ by an equivalent cocycle, we may therefore assume that the $Y_1$-modes in this range agree with the particular solution constructed in Section 4.1 of \cite{KQ-1st-def}. In particular, for every $j\leq k-1$, the mode formula in Section 4.1 (4) of \cite{KQ-1st-def} states that $e(-1)^{def} e(-1)^j \one $ may be computed by
\begin{align*}
     \frac 1 2 \sum_{\beta=1}^j e(-1)^{\beta-1} \sum_{\alpha\geq 0} \binom{-1}{\alpha}\left( (e(\alpha)e(-1)\one)^{def}_{-1+1-\alpha} + (e^{def}(\alpha)e(-1)\one)_{-1+1-\alpha}\right)e(-1)^{j-\beta}\one
\end{align*}
Since both $e(\alpha)e(-1)\one$ and $e^{def}(\alpha)e(-1)\one$ vanishes for every $\alpha \geq 0$, we conclude that for every $j\leq k-1$, 
\begin{align}
e(-1)^{def} e(-1)^j \one = 0. \label{e-def-e-i}
\end{align}
With (\ref{e-def-e-i}) and the induction hypothesis, (\ref{h0-def-e-i-1}) simplifies to 
\begin{align*}
    h^{def}(0) e(-1)^i\one
    = 0
\end{align*}
This proves the claim. 

\subsection{The terms with $e^{def}(-1)$} \label{e-1-def-sec}
We first show that (\ref{e-def-e-i}) also holds for $i=k$. Apply (\ref{h0-def-comm-condn}) to $e(-1)^k\one$. 
\begin{align*}
    & h^{def}(0) e(-1)^{k+1}\one - e(-1)h^{def}(0)e(-1)^k\one + h(0)e^{def}(-1)e(-1)^k\one - e^{def}(-1)h(0)e(-1)^k\one \nonumber\\
    = \ &  2 e^{def}(-1)e(-1)^k\one. 
\end{align*}
We analyze the left-hand-side. The first term vanishes because $e(-1)^{k+1}\one = 0$. The second term vanishes by the conclusion in Section \ref{h0-def-sec}. The fourth term is $2k e^{def}(-1)e(-1)^k\one$. Moving the fourth term to the right-hand-side results in 
$$h(0)e^{def}(-1)e(-1)^k\one = (2k+2)e^{def}(-1)e(-1)^k\one $$
So $e^{def}(-1)e(-1)^k\one$ is of conformal weight $k+1$ and of $h$-weight $2k+2$. Clearly, no such nonzero vector can exist. The conclusion then follows. In particular, the term $f(1)e^{def}(-1)e(-1)^{i-1}\one$ in (\ref{Recursion}) vanishes. 

Another term in (\ref{Recursion}) with $e^{def}(-1)$ is $e^{def}(-1)f(1)e(-1)^{i-1}$. We calculate that
\begin{align}
    f(1)e(-1)^{i-1}\one = \ & \sum_{j=0}^{i-2} e(-1)^j [f(1), e(-1)]e(-1)^{i-2-j}\one \nonumber \\
    = \ & \sum_{j=0}^{i-2} e(-1)^j (-h(0)+k) e(-1)^{i-2-j}\one \nonumber \\
    = \ & \sum_{j=0}^{i-2} e(-1)^j (-2(i-2-j)+k) e(-1)^{i-2-j}\one \nonumber \\
    = \ & ((-2i+4+k)(i-1) + (i-2)(i-1))e(-1)^{i-2}\one \nonumber \\
    = \ & (i-1)(k-i+2) e(-1)^{i-2}\one.
\end{align}
Therefore from (\ref{e-def-e-i}), we conclude that $e^{def}(-1)f(1)e(-1)^{i-1}\one = 0, 1\leq i\leq k+1$. 

\subsection{Proof of $c=0$} We are now ready to finish the proof. Since $e(-1)^{k+1}\one = 0$, with the results in Section \ref{h0-def-sec} and \ref{e-1-def-sec}, the recursion (\ref{Recursion}) with $i=k+1$ yields 
\begin{align*}
    0 = e(-1)f^{def}(1)e(-1)^{k}\one + c e(-1)^{k}\one 
\end{align*}
From the recursion (\ref{Recursion}) with $i=k$, we may obtain 
\begin{align*}
    0 = e(-1)^2f^{def}(1)e(-1)^{k-1}\one + 2c e(-1)^{k}\one 
\end{align*}
Repeatedly apply (\ref{Recursion}) with $i = k-1, ..., 1$, we obtain 
\begin{align*}
    0 = e(-1)^{k+1} f^{def}(1)\one + (k+1)\cdot c e(-1)^{k}\one 
\end{align*}
Since $f^{def}(1) \one =0$, we obtain that 
$$(k+1)\cdot c e(-1)^k \one = 0$$
Therefore, $c= 0$. We then conclude the following theorem. 
\begin{thm}
    A simple affine vertex operator algebra with positive integer level admits only trivial deformations. 
\end{thm}

\section{Affine $\mathfrak{sl}_2$ with admissible level $-4/3$}

We now study the case $V = L_{-4/3}(\mathfrak{sl}_2)$. 
It is known (see \cite{Adamovic} for details) that there is a singular vector of conformal weight 3 that is of the form 
$$s = -48 h(-1)e(-2)\one + 36 e(-1)^2 f(-1)\one - 6 h(-2)e(-1)\one + 9 h(-1)^2 e(-1)\one + 80 e(-3)\one$$
The relevant basis elements in $V^{-4/3}(\g)$ are 
$$h(-1)e(-2)\one, e(-1)^2f(-1)\one, h(-2)e(-1)\one, h(-1)^2e(-1)\one, e(-3)\one.$$
The strategy is basically the same as that in Section \ref{sec-int-level}, except that we need to consider both the $f^{def}(1) s$ and $h^{def}(1)s$ to obtain enough constraints for proving $c=0$. 

\subsection{$f^{def}(1)$-action on the singular vector}\label{f-def-modes} We note that the formulas in Section 4.1 (4) of \cite{KQ-1st-def} stay valid for every element of conformal weight 1 and 2. Using those formulas, together with the commutator condition as in Section 4.1, we compute that 
\begin{align*}
    f^{def}(1)h(-1)e(-2)\one = \ &  -f(1)h^{def}(-1)e(-2)\one.\\
    f^{def}(1)e(-1)^2 f(-1)\one = \ & 2c \cdot e(-1)f(-1)\one - f(1)e^{def}(-1)e(-1)f(-1)\one.\\
    f^{def}(1)h(-2)e(-1)\one= \ & c\cdot h(-2)\one-f(1)h^{def}(-2)e(-1)\one.\\
    f^{def}(1)h(-1)^2 e(-1)\one = \ & c\cdot h(-1)^2 \one-f(1)h^{def}(-1)h(-1)e(-1)\one.\\
    f^{def}(1)e(-3)\one = \ & 0. 
\end{align*}
Therefore, 
\begin{align*}
    f^{def} (1)s=  & -48f^{def} (1)h(-1)e(-2)\mathbf{1} +36f^{def} (1)e(-1)^{2} f(-1)\mathbf{1} -6f^{def} (1)h(-2)e(-1)\mathbf{1}\\
     & +9f^{def} (1)h(-1)^{2} e(-1)\mathbf{1} +80f^{def} (1)e(-3)\mathbf{1}\\
    =  & 48f(1)h^{def} (-1)e(-2)\mathbf{1} +72c\cdot e(-1)f(-1)\mathbf{1} -36f(1)e^{def} (-1)e(-1)f(-1)\mathbf{1} -6c\cdot h(-2)\mathbf{1}\\
     & +6\cdot f(1)h^{def} (-2)e(-1)\mathbf{1} +9\cdot c\cdot h(-1)^{2}\mathbf{1} -9\cdot f(1)h^{def} (-1)h(-1)e(-1)\mathbf{1}
\end{align*}

\subsection{$h^{def}(1)$-actions on the singular vector}\label{h-def-modes} With a similar process as in Section \ref{f-def-modes}, we compute that
\begin{align*}
    h^{def}(1)h(-1)e(-2)\one = \ & -h(1)h^{def}(-1)e(-2)\one+2c \cdot e(-2)\one.\\
    h^{def}(1)e(-1)e(-1)f(-1)\one = \ & -h(1)e^{def}(-1)e(-1)f(-1)\one\\
    h^{def}(1)h(-2)e(-1)\one= \ &  -h(1)h^{def}(-2)e(-1)\one. \\
    h^{def}(1)h(-1)h(-1)e(-1)\one = \ &  4c \cdot h(-1)e(-1)\one-h(1)h^{def}(-1)h(-1)e(-1)\one. \\
    h^{def}(1)e(-3)\one = \ & 0. 
\end{align*}
Therefore, 
\begin{align*}
    h^{def} (1)s=  & -48h^{def} (1)h(-1)e(-2)\mathbf{1} +36h^{def} (1)e(-1)^{2} f(-1)\mathbf{1} -6h^{def} (1)h(-2)e(-1)\mathbf{1}\\
     & +9h^{def} (1)h(-1)^{2} e(-1)\mathbf{1} +80h^{def} (1)e(-3)\mathbf{1}\\
    =  & 48h(1)h^{def} (-1)e(-2)\mathbf{1} -96c\cdot e(-2)\mathbf{1} -36h(1)e^{def} (-1)e(-1)f(-1)\mathbf{1}\\
     & +6h(1)h^{def} (-2)e(-1)\mathbf{1} +36c\cdot h(-1)e(-1)\mathbf{1} -9h(1)h^{def} (-1)h(-1)e(-1)\mathbf{1}.
\end{align*}

\subsection{Proof of $c=0$} Note that since $h^{def}(-1)e(-1)\one = 0$, by $D$-commutator formula, it is necessary that 
$$h^{def}(-2)e(-1)\one + h^{def}(-1)e(-2)\one = 0.$$
Therefore, it suffices to set $h^{def}(-1)e(-2)\one, e^{def}(-1)e(-1)f(-1)\one, h^{def}(-1)h(-1)e(-1)\one$ as generic elements of the conformal weight 3. Also, to determine $c$ it suffices to consider only those terms with Cartan weight 2. So we set  
\begin{align*}
    h^{def}(-1)e(-2)\one = \ & a_1 e(-3)\one + a_2 e(-2)h(-1)\one + a_3 e(-1)h(-2)\one \\ & + a_4 e(-1)h(-1)^2\one + a_5 e(-1)^2 f(-1)\one, \\
    h^{def}(-1)h(-1)e(-1)\one = \ & b_1 e(-3)\one + b_2 e(-2)h(-1)\one + b_3 e(-1)h(-2)\one \\ 
    & + b_4 e(-1)h(-1)^2\one + b_5 e(-1)^2 f(-1)\one.\\
    e^{def}(-1)e(-1)f(-1)\one= \ & c_1 e(-3)\one + c_2 e(-2)h(-1)\one + c_3 e(-1)h(-2)\one \\
    & + c_4  e(-1)h(-1)^2\one + c_5 e(-1)^2 f(-1)\one, \\
\end{align*}
Plug them into the formula of $f^{def}(1)s$ obtained in Section \ref{f-def-modes} and $h^{def}(1) s$ obtained in Section \ref{h-def-modes}, and set them to be zero, we obtain the following five equations concerning the structural constants:
\begin{align}
    0 = \ & 84 a_3 + 168 a_4 - 28 a_5 - 18 b_3 - 36 b_4 + 6 b_5 - 72 c_3 - 144 c_4 + 24 c_5 + 12 c  \label{eqn-l1}\\
0 = \ & -42 a_2 - 56 a_4 + 9 b_2 + 12 b_4 + 36 c_2 + 48 c_4 + 9 c \label{eqn-l2} \\
0 = \ & -42 a_1 - 56 a_3 + 9 b_1 + 12 b_3 + 36 c_1 + 48 c_3 - 6 c  \label{eqn-l3} \\
0 = \ & 84 a_2 - 560 a_4 + 168 a_5 - 18 b_2 + 120 b_4 - 36 b_5 - 72 c_2 + 
 480 c_4 - 144 c_5 + 36 c \label{eqn-l4}\\
 0 = \ & 84 a_1 - 280 a_2 - 168 a_3 + 784 a_4 - 336 a_5 - 18 b_1 + 60 b_2 + 36 b_3 - 
 168 b_4 + 72 b_5 \nonumber\\
 & - 72 c_1 + 240 c_2 + 144 c_3 - 672 c_4 + 288 c_5 - 96 c\label{eqn-l5}
\end{align}
We apply Gaussian elimination. The operation 
$$(\ref{eqn-l5})+2(\ref{eqn-l3})-\frac{70\cdot 4}{42} (\ref{eqn-l2}) + \frac{280}{84}(\ref{eqn-l1})$$
eliminates $a_1, a_2, a_3$ from (\ref{eqn-l5}) and results in 
\begin{align}
    \frac 4 3 (1288 a_4 - 322 a_5 - 276 b_4 + 69 b_5 - 96 c - 1104 c_4 + 276 c_5) = 0\label{eqn-nl5}
\end{align}
To eliminate $a_4$ in (\ref{eqn-nl5}), we first eliminate $a_2$ from (\ref{eqn-l4}) by the operation 
$$(\ref{eqn-l4})+2(\ref{eqn-l2}),$$
which results in 
\begin{align}
    -6 (112 a_4 - 28 a_5 - 24 b_4 + 6 b_5 - 9 c - 96 c_4 + 24 c_5) = 0 \label{eqn-nl4}
\end{align}
Then, the operation 
$$(\ref{eqn-nl5}) + \frac{1288 \cdot 4}{112\cdot 6 \cdot 3}(\ref{eqn-nl4})$$
results in 
$$10c = 0 \Rightarrow c = 0.$$
\begin{thm}
    The simple affine VOA $L_{-4/3}(sl_2)$ admits only trivial deformations.  
\end{thm}

\section*{Authors' statement}

\noindent The manuscript has no associated data. There exists no conflict of interest. We used ChatGPT 5.6-Pro to detect the gaps and mistakes in the previous version. The proof of Theorem \ref{thm-3-1} in the current version is contributed by ChatGPT 5.6-Pro.

\noindent {\small \sc Department of Mathematics, University of Denver, Denver, CO 80210, USA}

\noindent {\em E-mail address}: andrew.linshaw@du.edu

\noindent {\small \sc School of Mathematics (Zhuhai), Sun Yat-Sen University, Zhuhai, Guangdong, China}

\noindent {\em E-mail address}: qifei@mail.sysu.edu.cn | fei.qi.math.phys@gmail.com

\end{document}